\documentclass[leqno,12pt]{article}
\newcommand{\sca}{\ensuremath{\mathcal{A}}}
\newcommand{\scb}{\ensuremath{\mathcal{B}}}
\newcommand{\scd}{\ensuremath{\mathcal{D}}}
\newcommand{\sch}{\ensuremath{\mathcal{H}}}
\newcommand{\scx}{\ensuremath{\mathcal{X}}}
\newcommand{\dera}{\ensuremath{Der(\mathcal{A})}}
\newcommand{\aq}{\ensuremath{\mathcal{A}_Q}}

\begin{document}
\begin{center}
{\huge \textbf{Universality of quantum symplectic structure}}

\vspace{.2in} by TULSI DASS

\vspace{.3in} \textbf{Abstract}\end{center}

 Operating in the framework of `supmech'( a scheme of mechanics which
aims at providing a concrete setting for the axiomatization of
physics and of probability theory as required in Hilbert's sixth
problem; integrating noncommutative symplectic geometry and
noncommutative probability in an algebraic setting,
 it associates, with every `experimentally accessible' system, a
 symplectic algebra and operates essentially as noncommutative
 Hamiltonian mechanics with some extra sophistication in the
 treatment of states) it is shown that interaction between
 systems can be consistently described only if either (i) all
 system algebras are commutative or (ii) all system algebras are
 noncommutative and have a quantum symplectic structure
 characterized by a \emph{universal} Planck type
real-valued constant of the dimension of action.

\vspace{.2 in} \begin{center}
\[ \begin{array}{l}
Like \  it \ or \ not \\
If \ you \ are \ noncommutative \\
You \ have \ no \ option \\
But \ to \ be \ quantum. \end{array} \] \end{center}

\newpage  \begin{center} \textbf{Introduction} \end{center}

 Two (closely related) great intellectual challenges before
theoretical physicists are :

 \vspace{.12in} \noindent (i) Construction of the most economical and complete
description of nature (theory of `everything');

 \vspace{.12in} \noindent (ii) Solution of Hilbert's sixth problem [23]
(axiomatization of physics and probability theory).

For solving both these problems, two possible strategies are :

\vspace{.12in} \noindent (a) Solve (i), then brush up the formalism
and axiomatize so as to solve (ii).

\vspace{.12in} \noindent (b) Solve (i) in such a manner that (ii) is
automatically solved (essentially integrating the two problems).

The author's preference is for (b), mainly because, in this case,
relatively clearer thinking about (and contact with) fundamentals is
expected to prevail.

The adoption of (b) instead of (a) (which reflects the prevalent
attitude) implies a change in outlook and priorities. It puts
greater emphasis on the development of an `appropriate' formalism.
Without entering into a detailed discussion about the term
`appropriate', we shall take it to mean that the formalism should be
reasonably broad-based so as to cover \emph{all} systems in nature,
it should employ mathematics best suited for the development of the
adopted ideas and concepts and should be self consistent.

In the present era in physics, quantum theory is believed to be
applicable to all systems in nature. As far as experimental
predictions are concerned, it has been eminently successful. It is,
however, in need of a satisfactory formalism which should be in the
nature of its autonomous development (as opposed to the traditional
practice of \emph{quantizing} classical systems) and which should
provide for a satisfactory treatment of measurements on quantum
systems without introducing ad hoc assumptions like the von Neumann
reduction postulate.

The desired `appropriate' formalism must do justice to the basic
features of quantum mechanics (QM) : the noncommutative kinematics
of observables  and its intrinsically probabilistic nature as
reflected in the behavior of quantum states. The latter aspect,
traditionally referred to as `quantum probability' has been explored
in several versions [30], [26], [33], [27], [25], [35], [1], [31],
[28]. The one best suited to our needs is the one [28] based on
complex, associative, unital (i.e. having a unit element) and not
necessarily commutative
*-algebras (henceforth referred to as ALGEBRAs). In this version,
quantum probability may be referred to as noncommutative
probability. (Not. Since the term `noncommutative measure theory'
has been used for the algebraic development based on von Neumann
algebras presented in, for example, Connes' book [7], one might take
`noncommutative probability' to mean its `normalized' sub-domain; we
shall, however, reserve this term for the more general algebraic
version of Ref.[28].

The two (mutually related) noncommutative developments relating to
observables and states may be jointly referred to as the
`noncommutative culture' of QM.

Heisenberg's [22] idea --that kinematics underlying QM must be based
on a noncommutative algebra of observables - was incorporated into a
scheme of mechanics (called `matrix mechanics') by Born, Jordan,
Dirac and Heisenberg [4], [12], [5]. The proper geometrical
framework for the construction of the`quantum Poisson brackets' of
this mechanics is provided by noncommutative symplectic geometry
based on the derivation -based differential calculus developed by
Dubois-Violette and coworkers [16], [19], [17], [18], [14]; the
latter will be referred to as DVNCG.

Both, the noncommutative probability and DVNCG employ ALGEBRAs which
are, therefore, the natural domain for the development of the
`noncommutative culture' mentioned above. It makes perfect sense to
develop a coherent scheme of mechanics integrating noncommutative
symplectic geometry and noncommutative probability in the setting of
ALGEBRAs. Such a mechanics (called `supmech') has been developed by
the author. It has QM and classical Hamiltonian mechanics as special
subdisciplines and is projected as the appropriate framework for an
autonomous development of QM. The detailed development of this
mechanics will be presented elsewhere [11]. Here we shall restrict
ourselves to a reasonably self-contained presentation of a
development (within the domain of supmech) of some special
theoretical interest : a consistent description of interaction
between systems in the supmech framework is possible only if either

\vspace{.1in} \noindent  (i) all the system ALGEBRAs are
commutative, or

\vspace{.1in} \noindent (ii) all system ALGEBRAs are non-commutative
and have a quantum symplectic structure characterized by a
\emph{universal} real-valued constant of the dimension of action.

\vspace{.1in} \noindent  The formalism, therefore, has a natural
place for the Planck constant as a universal constant --- just as
special relativity has a natural place for a universal speed. In
fact, the situation in supmech is somewhat better because, whereas
in special relativity, the existence of a universal speed is
\emph{postulated}, in supmech the existence of a universal Planck
like constant is \emph{dictated/predicted} by the formalism.

\noindent Plan of the paper. In section 1, a brief account of DVNCG
is given which includes a discussion of its generalization [8]
involving algebraic pairs ($\mathcal{A},\mathcal{X}$) where
$\mathcal{A}$ is an ALGEBRA and $\mathcal{X}$ a Lie subalgebra of
Der($ \mathcal{A}$) and of the mappings [8], [9] induced on
derivations by the *-algebra isomorphisms (analogues of the
push-forward and pull-back mappings induced by diffeomorphisms on
vector fields and differential forms. In section 2, the
`noncommutative culture' of Hilbert space QM is expressed in
algebraic terms [to conform to the noncommutative geometry (NCG)-
based developments of the next section]. In section 3, an outline of
the supmech formalism is presented adequate for the treatment of
interacting systems in supmech in the next section. The last section
contains some concluding remarks.

\vspace{.15in} \noindent \textbf{Acknowledgements}. The author
thanks M. Dubois-Violette and M.J.W. Hall for their critical
comments on Ref.[8] and Ref.[10] respectively, to R. Sridharan and
V. Balaji for helpful discussions and to Chennai Mathematical
Institute and Indian Statistical Institute (Delhi Center) for
support and research facilities.

\vspace{.2in} \begin{center} \textbf{1. Derivation based
differential calculus} \end{center}

 \noindent 1.1 \emph{Noncommutative differential forms}. The
central object in DVNCG is an ALGEBRA $\mathcal{A}$; its elements
will be denoted as A,B,...and the identity element as I. The
*-operation (or involution) $* : \mathcal{A} \rightarrow \mathcal{A}
$ is an antilinear mapping satisfying the relations
\begin{eqnarray*} (AB)^* = B^*A^*, \ \ \  (A^*)^* = A, \ \ \ I^* = I.
\end{eqnarray*}
An element $A \in \mathcal{A}$ is called hermitian if $A^* = A$. The
center $Z(\mathcal{A})$ of $\mathcal{A}$ is the set of those
elements of $\mathcal{A}$ which commute with all elements of
$\mathcal{A}$.

A derivation of $\mathcal{A}$ is a linear map $X : \mathcal{A}
\rightarrow \mathcal{A}$ such that $ X(AB) = X(A)B + A X(B).$
Introducing the multiplication operator $\mu$ on \sca \ defined as
$\mu (A) B = AB$, the condition that X is a derivation may be
expressed as
\begin{eqnarray}
X \circ \mu (A) - \mu (A) \circ X = \mu (X(A)). \end{eqnarray} The
set $Der(\mathcal{A})$ of all derivations of $\mathcal{A}$ is a Lie
algebra with the Lie bracket $[X,Y] = X \circ Y - Y \circ X$. The
inner derivations $D_A$ defined by $D_AB = [A,B]$ satisfy the
relation
\begin{eqnarray*} [D_A , D_B] = D_{[A,B]} \end{eqnarray*}
and constitute a Lie subalgebra $IDer(\mathcal{A})$ of
$Der(\mathcal{A})$.

In DVNCG it is implicitly assumed that the ALGEBRAs being employed
have a reasonably rich supply of derivations so that various
constructions involving them have a nontrivial content.

An involution * on $Der(\mathcal{A})$ is defined by the relation
$X^*(A) = [X(A^*)]^*$. We have the (easily verifiable) relations
\begin{eqnarray*} [X,Y]^* = [X^*,Y^*], \ \ \ (D_A)^* = -D_{A^*}.
\end{eqnarray*}

By a \emph{differential calculus} on $\mathcal{A}$ one means a
formalism involving differential form like objects on $ \mathcal{A}$
with analogues of exterior product, exterior derivative and
involution defined on them. For noncommutative $\mathcal{A}$, it is
not unique; a systematic discussion of the variety of choices may be
found in Ref.[17]. In applications of NCG one makes a choice
according to convenience. In DVNCG (which is best suited for a
geometrical treatment of QM) one employs a derivation-based
differential calculus in which the spaces of differential p-forms
are (a subclass---to be specified later---of) Chevalley-Eilenberg
p-cochain spaces $C^p(Der(\mathcal{A}), \mathcal{A})$ [36]. Such a
p-cochain $\omega$ is, for $p\geq 1$, a multilinear map of
$[Der(\mathcal{A})]^p$ into $\mathcal{A}$ which is skew-symmetric :
\begin{eqnarray*}
\omega (X_{\sigma(1)},.., X_{\sigma(p)}) = \kappa_{\sigma} \omega
(X_1,..,X_p)
\end{eqnarray*}
where $\kappa_{\sigma}$ is the parity of the permutation $\sigma$;
we have $C^0(Der(\mathcal{A}, \mathcal{A}) = \mathcal{A}$.

An involution * on the cochains is defined by the relation
$\omega^*(X_1,..,X_p)= [\omega(X^*_1,..,X^*_p)]^*$; $\omega$ is said
to be real (imaginary) if $ \omega^* = \omega (-\omega)$.

The exterior product
\[ \wedge : C^p(Der(\mathcal{A}), \mathcal{A})
\times C^q(Der(\mathcal{A}), \mathcal{A}) \rightarrow
C^{p+q}(Der(\mathcal{A}, \mathcal{A}) \] is defined as in the
commutative case :
\begin{eqnarray}
(\alpha \wedge \beta)( X_1,.., X_{p+q}) & = & \frac{1}{p!q!}
\sum_{\sigma \in S_{p+q}} \kappa_{\sigma} \alpha(X_{\sigma(1)},..,
X_{\sigma(p)}). \nonumber \\
 &  \   & . \beta(X_{\sigma(p+1)},..,X_{\sigma(p+q)}).
\end{eqnarray}
With this product, the $N_0$-graded vector space
(where $N_0$ is the
set of non-negative integers)
\begin{eqnarray*} C(Der(\mathcal{A}, \mathcal{A})
= \bigoplus _{p \geq 0} C^p(Der(\mathcal{A}, \mathcal{A})
\end{eqnarray*}
becomes a graded complex algebra.

The Lie algebra $Der(\mathcal{A})$ acts on itself and on
$C(Der(\mathcal{A},\mathcal{A})$ through \emph{Lie
derivatives}. For
each $Y \in Der(\mathcal{A})$, one defines linear
mappings $ L_Y :
Der(\mathcal{A}) \rightarrow Der(\mathcal{A})$ and $L_Y :
C^p(Der(\mathcal{A}), \mathcal{A}) \rightarrow
C^p(Der(\mathcal{A}),\mathcal{A})$ such that the
following three
conditions hold :
\begin{eqnarray} L_Y (A) = Y(A) \ \ \textnormal{for all} \ \ A \in
\mathcal{A } \end{eqnarray} \begin{eqnarray} L_Y[X(A)] = (L_YX)(A) +
X[L_Y(A)] \end{eqnarray} \begin{eqnarray}
 L_Y [\omega(X_1,..,X_p)] & = &
(L_Y \omega)( X_1,..,X_p) \nonumber \\
 & \ & + \sum_{i=1}^p \omega(X_1,..,X_{i-1}, L_Y X_i,..,X_p).
\end{eqnarray} The first two conditions give
\begin{eqnarray} L_Y(X) = [Y,X] \end{eqnarray}
which, along with the third, gives
\begin{eqnarray} (L_Y \omega)(X_1,..,X_p) & = &
Y[\omega(X_1,..,X_p)] \nonumber \\
 & \ & - \sum_{i=1}^{p}\omega(X_1,..,X_{i-1},[Y,X_i],..,X_p).
\end{eqnarray}
Some important properties of the Lie derivative are, in obvious
notation, \begin{eqnarray} [L_X,L_Y] = L_{[X,Y]} \end{eqnarray}
\begin{eqnarray}
L_Y(\alpha \wedge \beta) = (L_Y \alpha) \wedge \beta + \alpha \wedge
(L_Y \beta).
\end{eqnarray}

For any $X \in \dera$, we define the interior product $ i_X :
C^p(\dera, \sca) \rightarrow C^{p-1}(\dera,\sca)$ (for $p \geq 1)$)
by
\begin{eqnarray} (i_X \omega)(X_1,..,X_{p-1}) =
\omega(X,X_1,..,X_{p-1}) \end{eqnarray} and $i_X(A) = 0$ for all $A
\in \sca$. The following relations involving the Lie derivative and
the interior product hold (here $\alpha$ is a p-form)
\begin{eqnarray} i_X \circ i_Y + i_Y \circ i_X = 0 \\
i_X (\alpha \wedge \beta) = (i_X \alpha) \wedge \beta + (-1)^p
\alpha \wedge (i_X \beta) \\
L_X \circ i_Y -i_Y \circ L_X = i_{[X,Y]}. \end{eqnarray}

The exterior derivative $ d : C^p(\dera, \sca) \rightarrow
C^{p+1}(\dera,\sca)$ is defined through the relation
\begin{eqnarray} (i_X \circ d + d \circ i_X) \omega = L_X \omega.
\end{eqnarray}
This equation determines the operation of d on cochains of various
degrees recursively. For p = 0, it takes the form
\begin{eqnarray} (dA)(X) = X(A). \end{eqnarray}
and, for general $ p \geq 0$,\begin{eqnarray}
 \ \ \ & \ & (d \alpha)(X_0,X_1, .., X_p)  \nonumber \\
 & = & \sum_{i=0}^p (-1)^i X_i[\alpha(X_0,..\hat{X}_i,..,X_p)]
  \nonumber \\
   & \ & + \sum_{0\leq i< j \leq p} (-1)^j
\alpha(X_0,.., X_{i-1}, [X_i,X_j], X_{i+1},.. \hat{X}_j,..,X_p)
\end{eqnarray}
where the hat indicates omission. The exterior derivative satisfies
the nilpotency condition $ d^2 =0$ and the relations
\begin{eqnarray} d \circ L_Y = L_Y \circ d \\
d (\alpha \wedge \beta) = d \alpha + \alpha \wedge (d \beta).
\end{eqnarray}
 The nilpotency of d implies that the pair $ (C(\dera, \sca),d)$
constitutes a cochain complex. We shall call a cochain $\alpha$
closed if $d \alpha =0$ and exact if $ \alpha = d \beta $ for some
cochain $\beta$.

Following Ref.[17], we consider the subset $ \Omega (\sca)$ of
$C(\dera, \sca)$ consisting of Z(\sca)-linear cochains which means
the cochains $\alpha$ satisfying the condition
\begin{eqnarray} \alpha (..,K X,..) = K \alpha (..,X,..)
\end{eqnarray}
for all $ X \in \dera \ $ and $K \in Z(\sca)$. This subset is closed
under the d-operation as can be easily easily verified using the
relation
\begin{eqnarray} [X, KY] = X(K)Y + K [X,Y] \end{eqnarray}
for all $X,Y \in \dera$ and $K \in Z(\sca)$. We shall reserve the
term `differential forms' for elements of $\Omega(\sca)$. We have
\begin{eqnarray*}
\Omega(\sca) = \bigoplus_{p\geq0} \Omega^p(\sca) \end{eqnarray*}
with $ \Omega^0(\sca) = \sca.$ Elements of $\Omega^p(\sca)$ will be
called differential p-forms.

\vspace{.15in} \noindent \emph{1.2  Induced mappings on derivations
and differential forms}

\vspace{.12in} A *-algebra isomorphism $ \Phi : \sca \rightarrow
\scb$ induces a mapping $ \Phi_* : \dera \rightarrow Der(\scb)$
given by
\begin{eqnarray} (\Phi_*X)(B) = \Phi (X[\Phi^{-1}(B)])
\end{eqnarray} for all $X \in \dera$ and $B \in \scb$. It is the
analogue (and a generalization) of the mapping induced by a
diffeomorphism on vector fields and satisfies the expected relations
(with $\Psi : \scb \rightarrow \mathcal{C})$
\begin{eqnarray} (\Psi \circ \Phi)_* = \Psi_* \circ \Phi_*; \ \
\Phi_* [X,Y] = [\Phi_*X, \Phi_*Y]. \end{eqnarray} It is easily seen
that $\Phi_*$ is a Lie-algebra isomorphism.

The *-isomorphism $\Phi$ also induces a mapping \begin{eqnarray*}
\Phi^* : C^p( Der(\scb), \scb) \rightarrow C^p(\dera, \sca)
\end{eqnarray*}  given, in obvious notation, by
\begin{eqnarray} (\Phi^* \omega)(X_1,..,X_p) = \Phi^{-1}
[\omega(\Phi_* X_1,.., \Phi_*X_p)]. \end{eqnarray} These mappings
are analogues (and generalizations) of the pull-back mappings on
traditional differential forms induced by diffeomorphisms. It is
easily seen that the mapping $\Phi_*$ preserves Z(\sca)-linear
combinations of derivations and that $ \Phi^*$ maps differential
forms onto differential forms. The following expected relations hold
:
\begin{eqnarray} (\Psi \circ \Phi)^* = \Phi^* \circ \Psi^* \\
\Phi^*(\alpha \wedge \beta) = (\Phi^* \alpha) \wedge (\Phi^* \beta)
\\
\Phi^*(d \alpha) = d (\Phi^* \alpha). \end{eqnarray}

Let $\Phi_t : \sca \rightarrow \sca $ be a one-parameter family of
transformations (i.e. ALGEBRA-automorphisms) given, for small t, by
\begin{eqnarray*}  \Phi_t (A) \simeq A + tg(A) \end{eqnarray*}
where g is some linear mapping of \sca \ into itself. The condition
\[ \Phi_t(AB) = \Phi_t(A) \Phi_t(B) \]
gives g(AB) = g(A)B + Ag(B)
implying that g(A) = Y(A) for some $Y \in \dera$; we call Y the
infinitesimal generator of the one-parameter family $\Phi_t$. It is
easily verified that the infinitesimal transformations of
derivations and of p-forms induced by $\Phi_t$ are given by the
respective Lie derivatives :
\begin{eqnarray} (\Phi_t)_* X \simeq X + t L_Y X \\
(\Phi_t)^* \omega \simeq \omega - t L_Y \omega. \end{eqnarray}

\vspace{.15in} \noindent \emph{1.3 Symplectic structures}

\vspace{.1in} A \emph{symplectic structure} on an ALGEBRA \sca \ is
defined as a differential 2-form $\omega$ (the \emph{symplectic
form}) which is (i) closed and (ii) non-degenerate in the sense
that, for every $ A \in \sca$, there is a unique derivation $Y_A$ in
\dera [the (globally) Hamiltonian derivation corresponding to A]
such that
\begin{eqnarray} i_{Y_{A}} \omega = - dA. \end{eqnarray}
The pair $(\sca, \omega)$ is called a \emph{symplectic algebra}.

A \emph{symplectic mapping} from a symplectic algebra $(\sca,
\alpha)$ to another one $(\scb, \beta)$ is an ALGEBRA-isomorphism
(i.e a *-algebra isomorphism mapping the unit element of \sca \ to
the unit element of \scb) such that $\Phi^* \beta = \alpha$. A
symplectic mapping from a symplectic algebra onto itself will be
called a \emph{canonical/symplectic transformation}. The symplectic
form and its exterior powers are invariant under canonical
transformations.

Given a symplectic algebra $(\sca, \omega)$, the \emph{Poisson
bracket} (PB) of two elements A and B of \sca \ is defined as
\begin{eqnarray} \{ A, B \} = \omega (Y_A,Y_B) = Y_A(B) = - Y_B(A).
\end{eqnarray}
It obeys the Leibnitz rule :
\begin{eqnarray} \{ A, BC \} = Y_A(BC) & = & Y_A(B)C + B Y_A(C)
 \nonumber \\
& = &  \{A, B \}C + B \{ A, C \}. \end{eqnarray} As in the classical
case [41], we also have the other two properties of PBs :

\vspace{.12in} \noindent (i) The Jacobi identity holds :
\begin{eqnarray} 0 & = & \frac{1}{2} (d \omega)(Y_A, Y_B, Y_C)
                         \nonumber \\
                   & = & \{A, \{B, C \} \} + \{B, \{C, A \} \}
                          + \{C, \{ A, B \} \};
\end{eqnarray} this makes (along with bilinearity and antisymmetry
of the PBs) the pair $(\sca, \{, \})$ a Lie algebra.

\vspace{.12in} \noindent (ii) The corespondence $ A \rightarrow Y_A$
is a Lie -algebra homomorphism from the above Lie algebra into \dera
:
\begin{eqnarray} [Y_A, Y_B] = Y_{ \{A, B \}}. \end{eqnarray}

An element A of \sca \ can act, via $Y_A$, as the infinitesimal
generator of a one-parameter family of canonical transformations.
The change in $ B\in \sca$ due to such an infinitesimal
transformation is
\begin{eqnarray}
\delta B = \epsilon Y_A(B) = \epsilon \{ A, B \}. \end{eqnarray}

\vspace{.15in} \noindent \emph{1.4 Canonical symplectic structure on
`special' ALGEBRAs}

\vspace{.12in} An ALGEBRA will be called \emph{special} if it has a
trivial center and if all its derivations are inner. The
differential 2-form $\omega_c$ defined on such an algebra \sca \ by
 \begin{eqnarray} \omega_c(D_A,D_B) = [A,B] \end{eqnarray}
 is said to be the \emph{canonical form} of \sca. (This differs from
 the definition in Ref.[16], [17] by a factor of i.) It is easily seen to
 be closed [the equation $(d\omega_c)(D_A, D_B, D_C) = 0 $ is
 nothing but the Jacobi identity for the commutator], imaginary (i.e.
 $\omega_c^* = - \omega_c) $ and dimensionless. For any $ A \in
 \sca$, the equation
 \begin{eqnarray*} \omega_c(Y_A, D_B) = - (dA)(D_B) = [A,B]
 \end{eqnarray*}  (for all $B \in \sca$) has the unique solution $Y_A =
 D_A$; this gives
 \begin{eqnarray} i_{D_A} \omega_c = -d A. \end{eqnarray}
The form $\omega_c$ defines, on \sca, the \emph{canonical symplectic
structure}; the corresponding PB is a commutator :
\begin{eqnarray} \{ A, B \} = D_A(B) = [A,B].
\end{eqnarray}

Using Equations (36) and (14), it is easily seen that the form
$\omega_c$ is \emph{invariant} in the sense that $L_X \omega_c = 0 $
for all $X \in \dera$. The invariant symplectic structure on the
algebra $M_n(C)$ of complex $n \times n$  matrices obtained  in Ref.
[19] is a special case of the canonical symplectic structure on
special ALGEBRAs described above.

If, on a special ALGEBRA \sca, instead of $\omega_c$, we take $
\omega = b \omega_c $ as the symplectic form (where b is a nonzero
complex number), we have
\begin{eqnarray} Y_A = b^{-1} D_A, \ \ \{ A, B \} = b^{-1}[A,B].
\end{eqnarray}
We shall see below that the so-called `quantum symplectic structure'
is such a symplectic structure with $ b = -i \hbar.$ Note that b
must be imaginary to make $\omega$ real. Just to have a convenient
name, we shall refer to the symplectic structure of the above sort
(for general non-zero b) as a quantum symplectic structure with
parameter b.

\vspace{.15in} \noindent \emph{ 1.5 A generalization of the
derivation-based differential calculus}

\vspace{.12in} A useful generalization of the formalism presented in
this section so far is obtained by restricting the derivations to a
Lie subalgebra $\scx$ of \dera; the central object in the whole
development will now be, instead of the ALGEBRA \sca, the pair
(\sca, \scx).  To get a feel for the implications of working with
such a pair, we consider a couple of examples, one `commutative' and
the other `noncommutative'.

\vspace{.12in} \noindent (i) $ \sca = C^{\infty}(R^3);$ \scx = the
Lie subalgebra of the Lie algebra $\scx (R^3)$ of smooth vector
fields on $R^3$ generated by the Lie differential operators $ L_j =
\epsilon_{jkl}x_k \partial_l$ for the SO(3)-action on $R^3$. These
differential operators act on the 2-dimensional spheres that
constitute the leaves of the foliation $R^3 -\{ (0,0,0) \} \cong S^2
\times R$. Employing the polar coordinates $(r,\theta, \phi)$ on
$R^3$ (which are obviously adapted to the above-mentioned
foliation), the variable r in the functions $f(r,\theta, \phi)$ in
$C^{\infty}(R^3)$ will remain unaffected by the derivations in \scx.
It follows that the restriction to the pair (\sca, \scx) in the
present case amounts to working on a leaf $(S^2)$ of the
above-mentioned foliation.

\vspace{.12in} \noindent (ii) $\sca = M_4(C)$, the algebra of
complex $4 \times 4$ matrices. The vector space $C^4$ on which these
matrices act serves as the carrier space of the spin $s= 3/2$
projective irreducible representation of the rotation group SO(3).
Denoting by $S_j (j = 1,2,3) $ the representatives of the generators
of the Lie algebra so(3) in this representation, let \scx \ be the
real Lie algebra generated by the inner derivations $D_{S_j}
(j=1,2,3)$. In the treatment of spin dynamics of a spin s = 3/2
object, one will effectively be using the pair (\sca, \scx).

\vspace{.1in}  In the generalized derivation-based differential
calculus based on a pair (\sca, \scx), one has the derivations
restricted to \scx \ and the p-cochains are those in the space
$C^p(\scx, \sca)$; the corresponding differential p-form space will
be denoted as $\Omega^p(\scx, \sca)$. Obviously $\Omega^p(\dera,
\sca) \equiv \Omega^p(\sca)$.

To define the induced mappings $\Phi_*$ and $\Phi^*$ in the present
context, one should employ a `pair isomorphism' $\Phi : (\sca,\scx)
\rightarrow (\scb, \mathcal{Y})$ which consists of an
ALGEBRA-isomorphism $\Phi : \sca \rightarrow \scb$ such that the
induced Lie algebra isomorphism $\Phi_* : \dera \rightarrow
Der(\scb)$ restricts to an isomorphism of \scx \ onto $\mathcal{Y}$.
Various properties of the induced mappings hold as before.

Given a one-parameter family of transformations (i.e.
pair-automorphisms) $ \Phi_t : (\sca, \scx) \rightarrow (\sca,
\scx)$, the condition $(\Phi_t)_* \scx \subset \scx$ implies that
the infinitesimal generator Y of $\Phi_t$ must satisfy the condition
$[Y,X] \in \scx$ for all $X\in \scx$. In practical applications, one
will generally have $Y \in \scx$ which obviously satisfies the
above-mentioned condition.

The concept of a symplectic algebra $(\sca, \omega)$ is now
generalized to that of a `generalized symplectic algebra' $(\sca,
\scx, \omega)$ where now $\omega \in \Omega^2(\scx, \sca)$. The
non-degeneracy condition on $\omega$ now demands, for a given $A \in
\sca$, the existence of a unique derivation $Y_A \in \scx$ such that
Eq.(29) holds. A symplectic mapping $\Phi : (\sca, \scx, \alpha)
\rightarrow (\scb, \mathcal{Y}, \beta)$ is now an
ALGEBRA-isomorphism $\Phi : \sca \rightarrow \scb$ such that the
induced mapping $\Phi_*$ restricts to an isomorphism of \scx \ onto
$\mathcal{Y}$ and $\Phi^* \beta = \alpha.$

\vspace{.2in} \begin{center} \textbf{2. The noncommutative culture
of quantum mechanics;  \\
the quantum symplectic structure}
\end{center}

 In this section, we shall present the traditional formalism of QM
in a not-so-familiar algebraic setting so as to obtain a useful
characterization of its `noncommutative culture'.

We start by considering the QM of a non-relativistic spinless
particle. The central object in it is the Hilbert space $\sch =
L^2(R^3,dx)$ of complex square-integrable functions on $R^3$. The
fundamental observables of such a particle are the Cartesian
components $X_j, P_j (j=1,2,3)$ of position and momentum vectors
which are self-adjoint linear operators represented, in the oft-used
Schr$\ddot{o}$dinger representation, as
\begin{eqnarray}
(X_j \phi)(x) = x_j \phi(x); \ \ (P_j \phi)(x) = -i \hbar
\frac{\partial \phi}{\partial x_j}. \end{eqnarray} These  operators
satisfy the canonical commutation relations (CCR)
\begin{eqnarray} [X_j, X_k] = 0 = [P_j, P_k]; \ \ \   [X_j, P_k] =
i \hbar I  \ \ (j,k = 1,2,3) \end{eqnarray} where I is the unit
operator. The functions $\phi$ in Eq.(39) must be restricted to a
suitable dense domain \scd \ in \sch \ which is generally taken to
be the space $\mathcal{S}(R^3)$ of Shwartz functions. Other
operators appearing in QM of the particle belong to the algebra \sca
\ generated by the operators $X_j, P_j $ (j= 1,2,3) and I [subject
to the CCR (40)]. The space $ \scd = \mathcal{S}(R^3)$ is clearly an
invariant domain for all elements of \sca. Defining a *-operation on
\scd \ by $ A^* = A^{\dagger}|\scd$, the Hermitian elements of \sca
\ represent the general observables of the particle.

A normalized element $\psi$ of \scd \ represents (up to a phase
factor)  a pure  state of the particle. Given the particle in this
state, the quantity
\begin{eqnarray} p(\Delta) \equiv \int_{\Delta} |\psi(x)|^2 dx
\end{eqnarray} (where $\Delta$ is a Borel set in $R^3$) is
interpreted as the probability that the particle lies in the domain
$\Delta$. For any observable $ A \in \sca$, the quantity
\begin{eqnarray}
<A>_{\psi} = (\psi, A \psi) \equiv \int \psi^*(x) (A \psi)(x) dx
\end{eqnarray} represents the expectation value of A in the state
$\psi$. With a suitable topology on the algebra \sca \ [15],the
quaantity $ \omega_{\psi} \equiv <.>_{\psi}$ of Eq.(42)  can be
considered as a continuous linear functional on \sca \ which is (i)
positive (which means  $\omega_{\psi}(B^*B) \geq 0 \ \ \forall B \in
\sca$) and (ii) normalized (i. e. $\omega_{\psi}(I) =1$). The set
$\mathcal{S}(\sca)$ of continuous positive linear functionals on
\sca \ is closed under convex combinations [i.e. $ \omega_i \in
\mathcal {S}(\sca) \Rightarrow \sum_i p_i \omega_i \in
\mathcal{S}(\sca)$ with $ p_i \geq 0, \sum_i p_i = 1$]. A nontrivial
convex combination of pure states is called a mixed state or
mixture.

It should now be easy to understand that a reasonably satisfactory
way of presenting the traditional formalism of QM of a system (which
permits free use of unbounded observables) is to associate, with a
quantum system S, a \emph{quantum triple} $(\sch, \scd, \aq)$ where
\sch \ is a complex, separable Hilbert space (which may or may not
be finite dimensional), \scd \ a dense linear domain in \sch (which
is obviously equal to \sch \ when \sch \ is finite dimensional) and
\aq \  an algebra of linear operators which, along with their
adjoints, have \scd \ as an invariant domain. For any $A \in \aq$,
we define its conjugate as $A^* = A^{\dagger}|\scd$ (thus defining
an involution * on \aq). Observables of the system are the Hermitian
elements of \aq.

For systems where a set of fundamental observables can be identified
(like the one considered above), the algebra \aq \ is the one
generated by the fundamental observables (and I) subject to
appropriate commutation relations.

States of S are those density operators $\rho$ such that
\begin{eqnarray}
<A>_{\rho} = Tr(\rho A) \end{eqnarray}
 is defined for all $A \in \aq$. For an observable A, the real
 quantity $<A>_{\rho}$ represents the expectation value of A when S
 is in the state $\rho$.  [Note. By states we strictly mean
 physical states so that expectation values of all observables are
 defined in all states.] Pure states are represented (up to  phase
 factors of modulus one) by normalized vectors $\psi \in \scd$ such
 that $<A>_{\psi} = (\psi, A \psi)$ The density operator
 corresponding to a state $\psi$ is $|\psi><\psi|$ in the Dirac
 notation.

 Dirac bra and ket spaces can be introduced in terms of Gelfand
 triples [20] based on the pair (\sch, \scd); we shall, however, skip the
 details.

 When the algebra \aq \ is `special' (in the sense defined in
 section 1), one has a canonical form $\omega_c$ defined on it
 [see Eq.(35)]. The \emph{quantum symplectic structure} is defined
 on \aq \ by employing the \emph{quantum symplectic form}
 \begin{eqnarray} \omega_Q = -i \hbar \omega_c. \end{eqnarray}
 Note that the factor i serves to make $\omega $ real and $\hbar$
 to give it the dimension of action (which is the correct dimension
 of a symplectic form in mechanics). The minus sign is a matter of
 convention. Eq.(38) now gives the \emph{ quantum Poisson bracket}
 \begin{eqnarray} \{ A, B \}_Q = (-i \hbar)^{-1} [A,B].
 \end{eqnarray}

 When the algebra \aq \ has both inner and outer derivations, one
 can employ the generalized symplectic algebra $( \aq, IDer(\aq),
 \omega_Q)$. Again, we have, for a given $A \in \aq$, $Y_A =
 (-i\hbar )^{-1} D_A$ and the quantum PB of Eq.(45).

 A nontrivial center in \aq \ indicates the presence of
 superselection rules and/or external fields. We shall skip details
 on these matters.

\vspace{.2in} \begin{center} \textbf{3. The formalism of supmech}
\end{center}

 As mentioned earlier, supmech is an algebraic scheme of mechanics
synthesizing noncommutative symplectic geometry and noncommutative
probability. Most developments in it are parallel to those in
classical Hamiltonian mechanics; in fact, it is essentially
noncommutative Hamiltonian statistical mechanics with some extra
sophistication in the treatment of states. In the detailed treatment
in Ref. [11], the basic system algebra is taken to be a superalgebra
(so as to provide a unified treatment of bosonic and fermionic
objects/entities); here, however, we shall restrict ourselves to the
simpler non-super version.

We shall call `experimentally accessible systems' those on which
repeatable experiments can be performed. For such systems, the
statistical analysis of experiments can be done with the traditional
frequency interpretation of probability. The universe as a whole and
large subsystems of it on a cosmological scale obviously do not
belong to this class. As of now, supmech has been developed only for
experimentally accessible systems.

The essential points in the development of supmech are listed below.

\vspace{.1in} \noindent \textbf{1.}\emph{The system algebra}.
Supmech associates with an experimentally accessible system S an
ALGEBRA \sca (its elements will be denoted as A,B,C,...). Hermitian
elements of \sca \ represent observables of S.We denote by
$\mathcal{O}(\sca)$ the set of all observables in \sca. (In fact,
\sca \ is assumed to be a locally convex algebra; we shall, however,
not treat the topological aspects here.)

\vspace{.1in} \noindent \textbf{2.}\emph{States}. States of \sca
(denoted by the letters $\phi, \psi,..$ )are defined as (continuous)
positive linear functionals which are normalized [i.e. $\phi(I) = 1$
where I is the unit element of \sca]. The set $\mathcal{S}(\sca)$ of
states of \sca \ is clearly closed under convex combinations
(weighted sums). Those states which cannot be represented as
nontrivial convex combinations are called pure. The set of pure
states of \sca \ is denoted as $\mathcal{S}_1(\sca)$. For any $A \in
\sca$ and $\phi \in \mathcal{S}(\sca)$, the quantity $\phi (A)$ is
to be interpreted as the expectation value of A in the state $\phi$.
When $A \in \mathcal{O}(\sca), \phi (A)$ is, of course, real.

\vspace{.1in} \noindent \textbf{3}. \emph{Compatible completeness of
observables and pure states}.  The pair
\[ (\mathcal{O}(\sca), \mathcal{S}_1(\sca)) \]
is assumed to be `compatibly complete' in the sense that \\
(i) given $ A,B \in \mathcal{O}(\sca), A \neq B$, there must be a
pure state $\phi$ such that $\phi(A) \neq \phi(B)$; \\
(ii) given two different pure states $\phi, \psi$, there must be an
observable A such that $\phi(A) \neq \psi(A)$.

\noindent  We shall refer to this as the \emph{CC condition}.

 \vspace{.1in} \noindent \textbf{4.}
\emph{Symplectic structure on the system algebra}. The system
algebra is assumed to have a symplectic structure provided by a
symplectic form $\omega$. Symmetries of the formalism (the analogues
of canonical transformations in classical Hamiltonian mechanics and
unitary transformations in QM) are canonical transformations of the
symplectic algebra $(\sca, \omega)$.

\vspace{.1in} \noindent  Note. The author has not opted for the
economy that could be obtained by combining items (\textbf{1}) and
(\textbf{4}) and introducing a system algebra directly as a
symplectic algebra because the first three items above constitute a
concrete unit serving a special purpose. [See remark (4) in the last
section.]

 \vspace{.12in} \noindent \textbf{5.} \emph{ Action of
canonical transformations on states}. Denoting the algebraic dual of
the algebra \sca \ by $\sca^*$, an automorphism $ \Phi : \sca
\rightarrow \sca$ induces the dual/transpose mapping $\tilde{\Phi} :
\sca^* \rightarrow \sca^*$ such that, in obvious notation,
\begin{eqnarray} \tilde{\Phi}(\phi)(A) = \phi(\Phi(A)) \ \
\textnormal{or} \ \ < \tilde{\Phi}(\phi), A> = < \phi, \Phi(A)>
\end{eqnarray}
where $ <,> $ denotes the dual space pairing. The mapping
$\tilde{\Phi}$ maps states (which form a subset of $\sca^*$) onto
states. To see this, note that \begin{eqnarray*}   &(i)& \ \
[\tilde{\Phi}(\phi)](A^*A) = \phi( \Phi(A^*A)) = \phi[\Phi(A)^*
\Phi(A)] \geq 0; \nonumber \\
& (ii) & \ \  [\tilde{\Phi}(\phi)](I) = \phi[\Phi(I)] = \phi(I)= 1.
\end{eqnarray*}

\noindent The linearity of $\tilde{\Phi}$ (on $\sca^*$) ensures that
it preserves convex combinations of states. In particular, it maps
pure states onto pure states. We have, therefore, a bijective
mapping $\tilde{\Phi} : \mathcal{S}_1(\sca) \rightarrow
\mathcal{S}_1(\sca)$.

When $\Phi$ is a canonical transformation, we have, for $X,Y \in
\dera$,
\begin{eqnarray*}
\omega(X,Y) = (\Phi^* \omega)(X,Y) = \Phi^{-1}[\omega(\Phi_*X,
\Phi_*Y)]
\end{eqnarray*}
giving
\begin{eqnarray}
\Phi[\omega(X,Y)] = \omega(\Phi_*X, \Phi_*Y).
\end{eqnarray}
Taking expectation value of each side in the state $\phi$, we get
\begin{eqnarray}
\tilde{\Phi}(\phi)[\omega(X,Y)] = \phi[\omega(\Phi_*X, \Phi_*Y)].
\end{eqnarray}
Defining $\omega_{\Phi}$ by
\begin{eqnarray} \omega_{\Phi}(X,Y) = \omega(\Phi_*x, \Phi_*Y)
\end{eqnarray} we can write Eq.(48) as
\begin{eqnarray} (\tilde{\Phi} \phi)[\omega(.,.)] =
\phi[\omega_{\Phi}(.,.)]. \end{eqnarray}

When $\Phi$ is an infinitesimal canonical transformation generated
by $G \in \sca$, we have
\begin{eqnarray} [\tilde{\Phi}(\phi)](A) \simeq \phi (A + \epsilon
 \{ G, A \}). \end{eqnarray} Putting $\tilde{\Phi}(\phi) = \phi +
\delta \phi$, we have
\begin{eqnarray} (\delta \phi)(A) = \epsilon \phi( \{G, A \}).
\end{eqnarray}

\vspace{.12in} \noindent \textbf{6.} \emph{Dynamics}. Dynamics of
the system is described by a one-parameter family $\Phi_t$ of
canonical transformations generated by an observable H called the
\emph{Hamiltonian}; the triple $ (\sca, \omega, H)$ will be called a
\emph{supmech Hamiltonian system}. As in QM or classical statistical
mechanics, there are two standard ways of describing dynamics
corresponding to the choice of making the evolution mappings act on
observables (Heisenberg type picture) or states
(Schr$\ddot{o}$dinger type picture); the two pictures are related as
[writing $\Phi_t(A) = A(t)$ and $\tilde{\Phi}_t(\phi) = \phi(t)$]
\begin{eqnarray} <\phi(t), A> = <\phi, A(t)>. \end{eqnarray}
In the Heisenberg type picture we have
\begin{eqnarray*} dA(t) = A(t+dt)-A(t) \simeq Y_H[A(t)] dt
\end{eqnarray*} giving the \emph{Hamilton's equation} of supmech :
\begin{eqnarray}
\frac{dA(t)}{dt} = Y_H[A(t)] = \{ H, A(t) \}. \end{eqnarray} In the
Schr$\ddot{o}$dinger type picture, Eq.(52) with $\Phi = \Phi_t$
gives the \emph{Liouville equation} of supmech :
\begin{eqnarray} \frac{d\phi(t)}{dt}(A) = \phi(t)(\{H, A \}) \ \
\textnormal{or}  \ \ \frac{d\phi(t)}{dt}(.) = \phi(t) (\{ H, . \}).
\end{eqnarray}

\vspace{.12in} \noindent \textbf{7.} \emph{Classical Hamiltonian
mechanics and QM as subdisciplines of supmech}.

\vspace{.1in} \noindent (i) \emph{Classical Hamiltonian mechanics}.
Traditionally developed in the framework of a symplectic manifold
$(M, \omega_{cl})$[41], it can be treated in supmech by taking $\sca
= C^{\infty}(M, C) \equiv \sca_{cl}$, the commutative algebra of
smooth complex-valued functions on the phase space M. The
observables of this systems are the subclass of real-valued
functions. For the algebra $\sca_{cl}$, the derivations  are the
smooth vector fields and the differential forms of section (1) are
the traditional differential forms on the manifold M. The symplectic
structure on $\sca_{cl}$ is given by the classical symplectic form
on M given, in standard notation, by
\[  \omega_{cl} = \sum_{j=1}^{n} dp_j \wedge dq^j \]
where dim (M) = 2n. Writing, in terms of the general local
coordinates $ \xi^a$ (a= 1,..,2n) on M, $ \omega_{cl} =
(\omega_{cl})_{ab} d\xi^a \wedge d\xi^b$,  the supmech Poisson
bracket on $\sca_{cl}$ is the classical Poisson bracket on M :
\begin{eqnarray} \{ f, g \}_{cl}  =  \omega_{cl}^{ab}
\frac{\partial f}{ \partial \xi^a} \frac{\partial g}{\partial \xi^b}
    =  \sum_j (\frac{\partial f} {\partial p_j}\frac{\partial
  g}{\partial q^j} - \frac{\partial f}{\partial q^j} \frac{\partial g}
  {\partial p_j}) \end{eqnarray}
where $(\omega_{cl}^{ab})$ is the inverse of the matrix
$((\omega_{cl})_{ab})$. The supmech Hamilton equation (54) is, in
the present context, the traditional Hamilton's equation
\begin{eqnarray} \frac{df}{dt} = \{ H_{cl}, f \}_{cl}.
\end{eqnarray}

States of $\sca_{cl}$ are probability measures on M; in obvious
notation, they are of the form $ \phi_{\mu}(f) = \int_M f d \mu.$
Pure states are Dirac measures (or, equivalently, points of M)
$\mu_{\xi_0}$ for which $\phi_{\xi_0}(f) = f(\xi_0)$.

The pair $(\mathcal{O}(\sca_{cl}), \mathcal{S}_1(\sca_{cl})) $ of
classical observables and pure states is easily sen to be compatibly
complete : Given two real-valued functions on M, there is a point of
M at which they take different values and, given two different
points of M, there is a real-valued function on M which takes
different values at those points.

In ordinary mechanics, only pure states are used. More general
states are used in classical statistical mechanics where, in most
applications, they are taken to be represented by densities on M [$d
\mu = \rho(\xi) d \xi$ where $d \xi = dq dp $ is the Liouville
volume element on M]. The state evolution equation of supmech gives,
in the present context,
\begin{eqnarray*} \int_M (\frac{\partial \rho(\xi,t)}{\partial
t}(\xi) f(\xi) d \xi = \int_M \rho (\xi,t) \{ H, f \}_{cl}(\xi) d
\xi. \end{eqnarray*} Taking $M = R^{2n}$, noting that the density
$\rho $ must vanish at infinity and performing a partial
integration, the right hand side becomes $\int_M \{ \rho,H \}_{cl} f
d \xi$ giving the traditional Liouville equation
\begin{eqnarray} \frac{\partial \rho} { \partial t} = \{\rho, H
\}_{cl}. \end{eqnarray}

\vspace{.12in} \noindent (ii) \emph{Quantum mechanics}. Most of the
needful has already been done in the previous section. Given a
quantum triple $(\sch, \scd, \sca_Q)$, the supmech system algebra is
to be taken as $\sca_Q$. The familly of pure states consists of unit
rays corresponding to vectors in \scd. The condition of compatible
completeness of the pair $(\sca_Q, \mathcal{S}_1 (\sca_Q))$ is
easily verified : \\
(i) Given $A,B \in \mathcal{O}(\sca_Q)$ and $ (\psi, A \psi) =
(\psi, B \psi)$ for all $\psi \in \scd$, we have $(\phi, A \psi) =
(\phi, B \psi)$ for all $\phi, \psi \in \scd$ implying A = B. [Hint:
Consider the given equality with state vectors $(\phi +
\psi)/\sqrt{2}$ and $ (\phi + i \psi)/\sqrt{2}$.] \\
(ii) Given normalized vectors $\phi, \psi \in \scd$, and $(\phi, A
\phi) = (\psi, A \psi)$ for all $A \in \mathcal{O}(\sca_Q)$, the
equality $\phi = \psi$ (up to a phase) can be seen by using the
given equality with A taken as the projection operators
corresponding to members of an orthonormal basis containg $\psi$ as
a member.

We have the quantum symplectic algebra $(\sca_Q, \omega_Q)$ and the
associated quantum Poisson brackets as in the previous section. The
supmech Hamilton equation (54) in the present case is clearly the
Heisenberg equation of motion
\begin{eqnarray*} \frac{dA(t)}{dt} = \{ H, A(t) \}_Q =
(-i\hbar)^{-1} [H, A(t)]. \end{eqnarray*} The supmech Liouville
equation (55) with the states given by density operators
$\omega_{\rho}(A) = Tr (\rho A)$ gives the `quantum Liouville
equation' (or the von Neumann equation)
\begin{eqnarray}
\frac{d \rho(t)}{dt} = (-i\hbar)^{-1}[\rho, H] = \{\rho(t), H \}_Q.
\end{eqnarray}

\vspace{.15in} \noindent \textbf{8.} \emph{Supmech as a framework
for an autonomous development of QM}

In the traditional development of QM, one generally \emph{quantizes}
classical systems. For example, to obtain the Schr$\ddot{o}$dinger
equation
\begin{eqnarray} i \hbar \frac{\partial \psi} {\partial t} =
[-\frac{\hbar^2}{2m} \nabla^2 +V] \psi \equiv H \psi \end{eqnarray}
in the traditional treatment of the QM of a nonrelativistic spinless
particle, one starts with the classical Hamiltonian
\begin{eqnarray} H = \frac {\mathbf{p}^2}{2m} + V, \end{eqnarray}
introduces the Hilbert space $\sch = L^2(R^3)$ of complex square
integrable functions, prescribes rules for the replacement of the
classical variables $x_j$ and $p_j$ by the operators $X_j$ and $P_j$
of Eq.(39) [thus obtaining the quantum Hamiltonian operator H of
Eq.(60)] and finally (taking clue from  the classical equation
$H_{cl} = E$ ), prescribes the rule for the evolution equation for
the Schr$\ddot{o}$dinger wave function $ \psi(\mathbf{x},t)$ in the
form $\hat{E} \psi = H \psi$ with $\hat{E} = i \hbar \frac{\partial}
{\partial t}$.

In Ref. [9], the need for an autonomous development of QM was
emphasized and some stringent criteria were laid down for such a
development. In the framework of supmech, it is possible to develop
the QM of particles autonomously satisfying those criteria [11]. We
give here an outline of the steps involved in the autonomous
development of the Schr$\ddot{o}$dinger equation (60). [The idea is
to define a particle as a localizable elementary system (which
involves a discussion of the action of the appropriate relativity
group on the system algebra and and of localizable systems), have a
systematic way to identify the fundamental observables of a particle
and obtain an expression for the Hamiltonian (infinitesimal
generator of time translations) in terms of the fundamental
observables (which can be done group theoretically [2]), and have a
systematic procedure to obtain a/the Hilbert space realization of
the relevant dynamics.]

\vspace{.12in} \noindent (i) One defines the \emph{Poisson action}
[41, [11] of a Lie group G on a symplectic algebra $(\sca, \omega)$
as an assignment, to every element $g \in G$, a canonical
transformation of the algebra such that the infinitesimal generators
(`hamiltonians') of one-parameter subgroups of the canonical
transformations have Poisson brackets in correspondence with the
commutation relations in the Lie algebra of G.

\vspace{.12in} \noindent (ii) The concept of a \emph{localizable
system} is introduced [as one which has a configuration space M (a
topological space) associated with it and it is meaningful to talk
about the probability of the system being localized in a Borel
subset of M] in which the concept of a position/configuration
observable naturally emerges. For systems with configuration space
$R^n$, the concept of \emph{concrete Euclidean-covariant
localization} is introduced in which one has the position
observables $X_j$ and the Euclidean group generators $P_j$ and
$M_{jk} (= -M_{kj})$ satisfying the standard Poisson bracket
relations.

\vspace{.12in} \noindent (iii) For the subclass of supmech systems
for which the concept of space and time and of a relativity scheme
are relevant, the appropriate relativity scheme is implemented
through the Poisson action of the corresponding relativity group
$G_0$ on the system algebra. In the nonrelativistic case (Galilean
relativity), the need for a Poisson action requires the replacement
of the Galilean group G by its projective group [2] $\hat{G}$ which
is a central extension of the universal covering group of G. The
additional generator corresponds to mass. In this manner, the
concept of mass appears naturally for the system at the fundamental
level.

For the implementation of a relativity scheme (with a relativity
group $G_0$), it is useful to introduce the concept of the
\emph{effective relativity group} $\hat{G}_0$ which is the universal
covering group $\tilde{G}_0$ of $G_0$ if the latter admits Poisson
actions and the projective group $\hat{G}_0$ if it does not.

\vspace{.12in} \noindent (iv) In supmech, an \emph{elementary
system} [for a given relativity scheme (or relativity group)] is
defined (generalizing and extending the treatments of elementary
systems by Wigner [40] and Alonso [2]) as a supmech triple $(\sca,
\omega, \mathcal{S}_1(\sca))$ such that  the effective relativity
group $\hat{G}_0$ has a Poisson action on the symplectic algebra
$(\sca, \omega)$ and a transitive action on the space
$\mathcal{S}_1(\sca)$ of pure states.

The fundamental observables of an elementary system are proposed to
be identified from the PBs of the `hamiltonians' coming from the
effective relativity group $\hat{G}_0$. For the Galilean elementary
systems, they turn out to be $M, X_j, P_j $ and $S_j$ (j=1,2,3)
corresponding, respectively, to mass, position, momentum and spin.
For a spinless particle they are $M, X_j$ and $P_j$. The observable
M (mass) has zero PBs with all other observables. It is, therefore,
a constant; its value m characterizes the elementary system and the
objects $X_j, P_j$ serve as kinematic observables. Simple group
theory leads to the following general expression (for massive
elementary systems)  for the generator H (the Hamiltonian ) of time
translations in terms of the fundamental observables :
\begin{eqnarray} H = \frac{\mathbf{P}^2}{2m} + V(\mathbf{X},
\mathbf{P}). \end{eqnarray}

\vspace{.12in} \noindent (v) A Hilbert space realization of the
supmech kinematics and dynamics of a system with noncommutative
algebra, if it exists, is very much desirable because, in such a
realization, the CC condition treated above is automatically
satisfied (as was seen in the subsection \textbf{7} above);
otherwise, one has to keep track of it separately. The existence of
a Hilbert space realization is, in fact, guaranteed by the CC
condition : there being a rich supply of (pure) states, one can
employ the GNS construction (the version of it best suited for us is
that of Ref.[24]) based on one of them to obtain a Hilbert space
representation of the algebra \sca. Such a representation is
generally not faithful; for example, if the state chosen is one with
zero expectation value for the kinetic energy (of a non-relativistic
particle), the momentum operator in the resulting Hilbert space
representation will be identically zero. The CC condition again
comes to the rescue; a faithful representation can be obtained by
taking an appropriate direct sum of the GNS representations of the
above sort.

All this trouble is, however, not necessary --- at least for a
system consisting of a single  particle. The condition of transitive
action of $\hat{G}_0$ on pure states implies that a Hilbert space
realization (in which pure states are vector states) must employ an
irreducible representation of this group. This, combined with the
points treated above, then ensures that the representation must be
the Schr$\ddot{o}$dinger representation. The probability
interpretation of Schr$\ddot{o}$dinger wave functions \emph{follows}
from the formalism. (This is because the essential relevant physics
is covered by the treatment of localizability above. This is, in
fact, very satisfying --- the probability interpretation of `$\psi$'
is no longer mysterious.)

The supmech evolution equation for pure states, with the Hamiltonian
of Eq.(62) (with V a function of X only in simple applications),
gives the traditional Schr$\ddot{o}$dinger equation. It should be
noted that the classical Hamiltonian or Lagrangian for a particle
was not used at any stage in this development.

\vspace{.12in} \noindent Note. Apart from ensuring an autonomous
development of QM and the interpretation of `$\psi$' above, a couple
of attractive features the formalism outlined above are :

\vspace{.12in} \noindent (1) The Planck constant $\hbar$ has to be
introduced `by hand' only once --- in the quantum symplectic form
(the most natural place to do it); its appearance at other
conventional places --- the canonical commutation relations (40),
the Heisenberg equation  and the Schr$\ddot{o}$dinger equation (60)
 --- is then automatic.

\vspace{.12in} \noindent (2) The Dirac bra-ket formalism (in its
rigorous version) appears naturally in the present setting. It is
this formalism --- and not von Neumann's formalism [30] employing
bounded observables --- which is used in most quantum mechanical
work.

\vspace{.2in} \begin{center} \textbf{4. Interacting systems in
supmech} \end{center}

 In this section, we shall consider, in the framework of supmech,
the interaction of two systems $S_1$ and $S_2$ described
individually as the supmech Hamiltonian systems $(\sca^{(i)},
\omega^{(i)}, H^{(i)})$ (i=1,2). We shall treat the coupled system
$S_1 + S_2$ also as a supmech Hamiltonian system. To this end,  we
associate, with the coupled system $S_1 + S_2$ the (algebraic)
tensor product algebra $ \sca = \sca^{(1)} \otimes \sca^{(2)}$. The
most important job in the present section is, given the symplectic
forms $ \omega^{(1)}$ and $\omega^{(2)}$ on $\sca^{(1)}$ and
$\sca^{(2)}$,  to determine the symplectic form and the PB on \sca.

\vspace{.12in} \noindent \emph{4.1 The symplectic form and PB on the
algebra $\sca = \sca^{(1)} \otimes \sca^{(2)}$}

\vspace{.12in} The algebra $\sca^{(1)}$ (resp. $\sca^{(2)})$ has, in
\sca, an isomorphic copy consisting of the elements ($A \otimes I_2,
A \in \sca^{(1)}$) (resp. $I_1 \otimes B, B\in \sca^{(2)}$) to be
denoted as $\tilde{\sca}^{(1)}$ (resp. $\tilde{\sca}^{(2)}$) where
$I_1$ and $I_2$ are the unit elements of $\sca^{(1)}$ and
$\sca^{(2)}$ respectively. We shall also use the notations
$\tilde{A}^{(1)} \equiv A \otimes I_2$ and $\tilde{B}^{(2)} \equiv
I_1 \otimes B$.

Objects in $\sca^{(i)}$ and $\tilde{\sca}^{(i)}$ are related through
the induced mappings corresponding to the isomorphisms $\Xi^{(i)} :
\sca^{(i)} \rightarrow \tilde{\sca}^{(i)}$ (i= 1,2) given by
$\Xi^{(1)}(A) = A \otimes I_2$ and $\Xi^{(2)}(B) = I_1 \otimes B$.
In particular \\
(i) The induced mapping $\Xi^{(1)}_* : Der (\sca^{(1)}) \rightarrow
Der (\tilde{\sca}^{(1)})$ gives $\Xi^{(1)} _*(X) = \tilde{X}^{(1)}$
where
\begin{eqnarray*}
\tilde{X}^{(1)} (\tilde{A}^{(1)}) = \Xi^{(1)}[X(A)] = X(A) \otimes
I_2.
\end{eqnarray*} Similarly, corresponding to $Y \in Der(\sca^{(2)})$,
we have $\tilde{Y}^{(2)} \in Der(\tilde{\sca}^{(2)})$ given by
$\tilde{Y}^{(2)}(\tilde{B}^{(2)}) = I_1 \otimes Y(B)$. \\
(ii) The induced mappings on 1-forms give, corresponding to the
1-forms $\alpha \in \Omega^1(\sca^{(1)})$ and $\beta \in
\Omega^1(\sca^{(2)})$, we have $\tilde{\alpha}^{(1)} \in
\Omega^1(\tilde{\sca}^{(1)})$ and $\tilde{\beta}^{(2)} \in
\Omega^1(\tilde{\sca}^{(2)})$ given by
\begin{eqnarray*} \tilde{\alpha}^{(1)}(\tilde{X}^{(1)}) =
\Xi^{(1)}[\alpha(([\Xi^{(1)}]^{-1})_* \tilde{X}^{(1)})] =
\Xi^{(1)}[\alpha(X)] = \alpha(X) \otimes I_2 \end{eqnarray*} and
$\tilde{\beta}^{(2)}(\tilde{Y}^{(2)}) = I_1 \otimes \beta(Y)$.
Similar formulas hold for the higher forms.

To obtain the general differential forms and the exterior derivative
on \sca, the most straightforward procedure is to obtain the graded
differential space $ (\Omega (\sca), d)$ as the tensor product [21]
of the graded differential spaces $(\Omega (\sca^{(1)}), d_1)$ and $
(\Omega(\sca^{(2)}, d_2).$ A differential k-form on \sca \ is of the
form (in obvious notation)
\begin{eqnarray*} \alpha_k = \sum_{i +j = k} \alpha_i^{(1)} \otimes
\alpha_j^{(2)}. \end{eqnarray*} The exterior derivative d on $\Omega
(\sca)$ is given by [here $\alpha \in \Omega^p(\sca^{(1)})$ and
$\beta \in \Omega(\sca^{(2)})$] \begin{eqnarray} d(\alpha \otimes
\beta) = (d_1 \alpha) \otimes \beta + (-1)^p  \alpha \otimes (d_2
\beta). \end{eqnarray}

Given the symplectic forms $\omega^{(i)}$ on $\sca^{(i)}$ (i= 1,2)
and stipulating that the symplectic form $\omega$ on \sca \ should
not depend on anything other than the objects $ \omega^{(i)}$ and
$I_{(i)}$ (i=1,2) (the `naturality'/`canonicality' assumption), the
only possible choice of $\omega$ is
\begin{eqnarray} \omega = \omega^{(1)} \otimes I_2 + I_1 \otimes
\omega^{(2)}. \end{eqnarray} To show that it is, indeed, a
symplectic form, we must show that it is (i) closed and (ii)
non-degenerate. Eq.(63) gives
\begin{eqnarray*} d \omega = (d_1 \omega^{(1)}) \otimes I_2 +
\omega^{(1)} \otimes d_2 (I_2) + d_1(I_1) \otimes \omega^{(2)} + I_1
\otimes d_2 \omega^{(2)} = 0 \end{eqnarray*} showing that $\omega $
is closed.

To show the non-degeneracy of $\omega$, we must show that, given $A
\otimes B \in \sca$, there exists a unique derivation $ Y = Y_{A
\otimes B}$ in \dera \ such that \begin{eqnarray} i_Y \omega = -d (A
\otimes B) & = & -(d_1A) \otimes B - A \otimes (d_2 B) \nonumber \\
           & = & i_{Y_A^{(1)}} \omega^{(1)} \otimes B + A \otimes
           i_{Y_B^{(2)}} \omega^{(2)}. \end{eqnarray}
The structure of Eq.(65) suggests that Y must be of the form
\begin{eqnarray} Y = Y_A^{(1)} \otimes \Psi_B^{(2)} + \Psi_A^{(1)}
\otimes Y_B^{(2)} \end{eqnarray} where $\Psi_A^{(1)}$ and
$\Psi_B^{(2)}$ are linear mappings on $\sca^{(1)}$ and $\sca^{(2)}$
respectively such that $\Psi_A^{(1)} (I_1) = A$ and $ \Psi_B^{(2)}
(I_2) = B.$ A general object of the form (66), however, need not be
a derivation of \sca; we must, therefore, impose the condition that
Y must be a derivation. Recalling Eq.(1) and denoting the
multiplication operators in $\sca^{(1)}, \sca^{(2)}$ and \sca \ by
$\mu_1, \mu_2$ and $\mu$ respectively, we have \begin{eqnarray} Y
\circ \mu (C \otimes D) - \mu(C \otimes D) \circ Y = \mu(Y(C \otimes
D)). \end{eqnarray} Noting that $\mu(C \otimes D) = \mu_1(C) \otimes
\mu_2(D)$, Eq.(67) with Y of Eq.(66) gives \begin{eqnarray} (
Y_A^{(1)} \circ \mu_1(C) ) \otimes ( \Psi_B^{(2)} \circ \mu_2(D) ) +
( \Psi_A^{(1)} \circ \mu_1(C) ) \otimes ( Y_B^{(2)} \circ
\mu_2(D) ) \nonumber \\
- (\mu_1(C) \circ Y_A^{(1)} ) \otimes ( \mu_2(D) \circ \Psi_B^{(2)})
- ( \mu_1(C) \circ \Psi_A^{(1)} ) \otimes (
\mu_2(D) \circ Y_B^{(2)} ) \nonumber \\
= \mu [ Y_A^{(1)} (C) \otimes \Psi_B^{(2)}(D) + \Psi_A^{(1)}(C)
\otimes Y_B^{(2)}(D)]. \end{eqnarray} Since $Y_A^{(1)}$ and
$Y_B^{(2)}$ are derivations, we must have \begin{eqnarray} Y_A^{(1)}
\circ \mu_1(C) - \mu_1(C) \circ Y_A^{(1)} = \mu_1 (Y_A^{(1)}(C)) =
\mu_1 (\{ A,C \}_1) \nonumber \\
Y_B^{(2)} \circ \mu_2(D) - \mu_2(D) \circ Y_B^{(2)} =
\mu_2(Y_B^{(2)}(D)) = \mu_2 ( \{ B,D \}_2). \end{eqnarray}

Putting $D = I_2$ in Eq.(68), we have [noting that $\mu_2(D) =
\mu_2(I_2) = id_2$, the identity mapping on $\sca^{(2)}$ and
$Y_B^{(2)}(I_2) = 0$]\begin{eqnarray*}  \ \ \ & \ & ( Y_A^{(1)}
\circ \mu_1(C) ) \otimes
\Psi_B^{(2)} + ( \Psi_A^{(1)} \circ \mu_1(C) ) \otimes Y_B^{(2)} \\
 & \ & - ( \mu_1(C) \circ Y_A^{(1)} ) \otimes \Psi_B^{(2)} - (
\mu_1(C)_ \circ Y_A^{(1)} ) \otimes Y_B^{(2)}  \\
 & = & \mu [Y_A^{(1)}(C) \otimes B] = \mu_1 (\{ A,C \}_1) \otimes
\mu_2(B) \end{eqnarray*}
 which, along with Eq.(69), gives
\begin{eqnarray} \ \ \ \ \ \mu_1(\{A,C \}_1) \otimes [\Psi_B^{(2)} - \mu_2(B)]
= [\mu_1(C) \circ \Psi_A^{(1)} - \Psi_A^{(1)} \circ \mu_1(C)]
\otimes Y_B^{(2)}.
\end{eqnarray}

Similarly, putting $C = I_1$ in Eq.(68), we get \begin{eqnarray} \ \
\ \  [\Psi_A^{(1)} - \mu_1(A) ] \otimes \mu_2(\{ B,D \}) = Y_A^{(1)}
\otimes [\mu_2(D) \circ \Psi_B^{(2)} - \Psi_B^{(2)} \circ \mu_2(D)].
\end{eqnarray}

Now, equations (71) and (70) give \begin{eqnarray} \Psi_A^{(1)} -
\mu_1(A) = \lambda_1 Y_A^{(1)}  \\
\mu_2(D) \circ \Psi_B^{(2)} - \Psi_B^{(2)} \circ \mu_2(D) =
\lambda_1 \mu_2( \{ B,D \}_2)  \\
\Psi_B^{(2)} - \mu_2(B) = \lambda_2 Y_B^{(2)}  \\
\mu_1(C) \circ \Psi_A^{(1)} - \Psi_A^{(1)} \circ \mu_1(C) =
\lambda_2 \mu_1( \{A,C \}_1) \end{eqnarray} where $\lambda_1$ and
$\lambda_2$ are complex numbers.

Equations (66), (72) and (74) give
\begin{eqnarray} Y & = & Y_A^{(1)} \otimes [\mu_2 (B) + \lambda_2
Y_B^{(2)}] + [\mu_1(A) + \lambda_1 Y_A^{(1)} ] \otimes Y_B^{(2)}
\nonumber \\
& = & Y_A^{(1)} \otimes \mu_2(B) + \mu_1(A) \otimes Y_B^{(2)} +
(\lambda_1 + \lambda_2) Y_A^{(1)} \otimes Y_B^{(2)}.
\end{eqnarray} Note that only the combination $(\lambda_1 +
\lambda_2) \equiv \lambda $ appears in Eq.(76). To have a unique Y,
we must obtain an equation fixing  $\lambda$ in terms of given
quantities.

Substituting for $\Psi_A^{(1)}$ and $\Psi_B^{(2)}$ from equations
(72) and (74) into equations (73) and (75) and using equations (69),
we obtain the equations
\begin{eqnarray} \lambda \mu_1 (\{ A,C \}_1) = \mu_1([C,A]) \ \
\textnormal{for all} \ \ A,C \in \sca^{(1)}  \\
\lambda \mu_2( \{ B,D \}_2) = \mu_2 ([D,B]) \ \ \textnormal{for all}
\ \ B,D \in \sca^{(2)}. \end{eqnarray} We have not one but two
equations of the type we have been looking for. This is a signal for
the emergence of nontrivial conditions (for the desired symplectic
structure on the tensor product algebra to exist).

 Let us consider the equations (77,78) for the various possible
 situations:

 \vspace{.12in} \noindent (i) Let the algebra $\sca^{(1)}$ be commutative.
 Assuming the PB $\{ , \}_1$ is nontrivial, Eq.(77) implies that
 $\lambda =0$. Then Eq.(78) implies that the algebra $\sca^{(2)}$ is
 also commutative. It follows that \\
(a) when both the algebras $\sca^{(1)}$ and $ \sca^{(2)}$ are
commutative, the unique Y is given by Eq.(76) with $\lambda =0$; \\
 (b) a `natural'/`canonical' symplectic structure  does not exist
 on the tensor product of a commutative and a noncommutative
 algebra.

\vspace{.12in} \noindent (ii) Let the algebra $\sca^{(1)}$ be
noncommutative. Eq.(77) then implies $\lambda \neq 0$ which, in
turn, implies, through Eq.(78), that the algebra $\sca^{(2)}$ is
also non-commutative [which is also expected from (b) above].
Equations (77, 78) now give \begin{eqnarray} \{ A, C \}_1 = -
\lambda^{-1} [A,C], \ \ \ \{ B, D \}_2 = - \lambda^{-1} [B,D]
\end{eqnarray} which shows that when both the algebras $\sca^{(1)}$
and $\sca^{(2)}$ are noncommutative, a `natural'/`canonical'
symplectic structure on their tensor product  exists if and only if
each algebra has a quantum symplectic structure with the \emph{same}
parameter (-$ \lambda $) , i.e.
\begin{eqnarray} \omega^{(1)} = - \lambda \omega^{(1)}_c , \ \ \
\omega^{(2)} = - \lambda \omega^{(2)}_c \end{eqnarray} where
$\omega^{(1)}_c$ and $\omega^{(2)}_c$ are the canonical symplectic
forms on the two algebras. It follows that \emph{all noncommutative
system algebras must have a universal quantum symplectic structure}.
Comparison of Eq.(80) with the quantum symplectic form (44) shows
that $ \lambda = i \hbar$.

In all the permitted cases, the PB on the algebra $ \sca =
\sca^{(1)} \otimes \sca^{(2)}$  is given by
\begin{eqnarray}  \{ A \otimes B, C \otimes D \} & = & \{ A, C \}_1
\otimes BD + AC \otimes \{ B,D \}_2 \nonumber \\
 & \ & + \lambda \{ A,C \}_1 \otimes \{ B,D \}_2 \end{eqnarray}
 where the parameter $\lambda$ vanishes in the commutative case; in the
noncommutative case, it is the universal parameter appearing in the
symplectic forms (80).

In Ref. [10], the following PB was reported for the tensor product
algebra \sca \ :
\begin{eqnarray} \ \ \ \ \{ A \otimes B, C \otimes D \} = \{ A,B \}_1
\otimes \frac{CD + DC}{2} + \frac{AC + CA}{2} \otimes \{ C,D \}_2.
\end{eqnarray} When both the algebras $\sca^{(1)}$ and $\sca^{(2)}$
are commutative, the equations (81) (with $\lambda = 0$) and (82)
are clearly the same. In fact, the same is also true when both the
algebras are noncommutative. To see this, it is adequate to note
that, using Eq.(79), we have \begin{eqnarray*}
 \lambda \{ A,C \}_1
\otimes \{ B,D \}_2  & = & [C,A] \otimes \{ B,D \}_2 = \{ A,C \}_1
                            \otimes [D,B] \nonumber \\
                     & = & \frac{CA - AC}{2} \otimes \{ B,D \}_2 +
                     \{ A,C
\}_1 \otimes \frac{DB - BD}{2}. \end{eqnarray*}

In Ref. [10], the PB of Eq.(82) was meant to be true for the general
case which includes the case when one of the two algebras is
commutative and the other noncommutative ( the mixed case). Shortly
after the paper in Ref.[10] appeared in the arXiv, M.J.W. Hall, in a
private communication to the author, pointed out that the `Poisson
bracket' (82) does not satisfy the Jacobi identity in some cases (as
shown, for example, in Ref.[6]). The present work is an outcome of
the efforts to clarify the situation in this matter.

The example of violation of the Jacobi identity belonged to the
mixed case. We now know that, in this case, Eq.(64) does not
represent a valid symplectic structure. The mistake in the earlier
work of the author consisted in not ensuring that the Y of Eq.(66)
is a derivation.

Comment on a possible generalized symplectic structure (of the type
mentioned in section 1.5) is being postponed to the last section
(item 5 there).

\vspace{.15in} \noindent \emph{4.2 Dynamics of the interacting
system}.  Given that the two systems $S_1$ and $S_2$ are represented
as supmech Hamiltonian systems as mentioned above, the coupled
system $ S_1 + S_2$ is also to be represented as a supmech
Hamiltonian system $ (\sca, \omega, H)$ with the symplectic form
$\omega$ as in Eq.(64) and the Hamiltonian given by
\begin{eqnarray} H = H^{(1)} \otimes I_2 + I_1 \otimes H^{(2)} +
H_{int}. \end{eqnarray} In most applications, the interaction
hamiltonian is of the form \begin{eqnarray} H_{int} = \sum_{i=
1}^{n} F_i \otimes G_i \end{eqnarray} where $F_i$ and $G_i$ are
observables of the two systems and the coupling constants have been
absorbed in these observables.

In the Heisenberg type picture, evolution of a typical observable
$A(t) \otimes B(t)$ is governed by the supmech Hamilton equation
\begin{eqnarray} \frac{d}{dt} [A(t) \otimes B(t)] & = & \{ H, A(t)
\otimes B(t) \} \nonumber \\
& = & \{ H^{(1)}, A(t) \}_1 \otimes B(t) + A(t) \otimes \{ H^{(2)},
B(t) \}_2  \nonumber \\
 & \ & + \{ H_{int}, A(t) \otimes B(t) \}. \end{eqnarray}

In Ref.[10], this formalism was applied to the treatment of
measurements in quantum mechanics taking $S_1$ to be the measured
quantum system and $S_2$ the apparatus (assumed macroscopic) treated
as a classical system and using the PB of Eq.(82). Since a
quantum-classical interaction is now not permitted, we must go back
to the original von Neumann idea [30] to treat the apparatus as a
quantum mechanical system. We are, however, not constrained to adopt
the von Neumann procedure [30], [38] of introducing vector states
for the pointer positions (which is the basic cause of all the
problems in quantum measurement theory). Here supmech offers an
advantage not available in  von Neumann's treatment. Since both
quantum and classical systems can be accommodated in the supmech
formalism, one can exploit the fact that the apparatus can be
described classically to a very good approximation. The best way to
do this is to use the phase space description of the QM of the
apparatus (the Weyl-Wigner-Moyal formalism [37], [39], [29])and then
go to the classical approximation (doing it all within the supmech
formalism). With such a modification, the program of Ref[10] goes
through successfully, justifying the final results obtained there.
We shall, however, skip the details here.

\vspace{.2in} \begin{center} \textbf{5. Concluding remarks}
\end{center}

 \noindent 1. Supmech permits two kinds of `worlds' : the
\emph{commutative world} in which all system algebras are
commutative and the \emph{noncommutative world} in which they are
all noncommutative. There is no restriction (as far as the
consistency of the supmech formalism is concerned) on the possible
symplectic structures on system algebras in the commutative world;
however, the system algebras in the noncommutative world must all
have a universal quantum symplectic structure. Since QM is known to
describe systems in nature substantially correctly, the real world
is, of course, noncommutative; systems in the commutative world can
appear only as approximations to those in the real quantum world.

\vspace{.12in} \noindent 2. The existence of a natural place for a
universal Planck-like constant in the formalism is an important
feature of supmech and deserves further comment and elaboration.

In physics, out of the three fundamental constants G (Newton's
constant of gravitation), c (the speed of light in vacuum) and
$\hbar$ (the Planck constant), the first (G) appears in the
statement of a universal law of nature (in Newton's law of
gravitation and in Einstein's gravitational field equation in the
general theory of relativity); the second (c) appears in classical
electromagnetic theory as the speed of electromagnetic waves in
vacuum, it is \emph{postulated} as a universal speed in special
relativity and maintains such existence in general relativity
through the equivalence principle. The last one ($\hbar$) was
introduced in the relation $ E = h \nu = \hbar \omega$ as the
proportionality constant between energy and frequency in the
hypothesized fundamental unit (`quantum') of energy in the energy
exchange between interacting systems; in the traditional development
of QM, it is put `by hand' in various equations --- the canonical
commutation relations, the Heisenberg's equation of motion and  the
Schr$\ddot{o}$dinger equation. As  has been already mentioned, QM is
need of a proper formalism. The fact that supmech, apart from its
geometrical setting and other appealing features, predicts the
existence of a universal Planck-like constant, is a strong
indication that the `right' formalism has been chosen for an
autonomous development of QM.

\vspace{.12in} \noindent 3. If one could construct a formalism in
which there are \emph{similar} natural places for three independent
dimensional parameters [say, $\hbar$, c  and \emph{l} (a fundamental
length)], it would constitute substantial progress towards
construction of the `theory of everything'. For this, one might try
to find sub-theories of supmech with natural places for c and
\emph{l} or, more generally, supmech-like theories which have the
above-mentioned feature of supmech for all the three parameters.

Emphasis on the word `similar' in the previous para means that the
other two universal constants should also appear as proportionality
constants in the choices of appropriate geometrical objects as
multiples of the corresponding `canonical' objects (recall $
\omega_Q = -i \hbar \omega_c$) --- or through some similar
compelling geometrical reasoning. If we relax this requirement, one
can find other ways of having reasonably `natural' looking places
for universal constants which may not have as profound implications
as the appearance of the parameter $\hbar$ had in supmech. For
example, one may choose to work in a spatial lattice of fundamental
spacing `a' and employ discrete evolution with step length `b' of
the evolution parameter (`time'); one can now take \emph{l} = a and
c = a/b. While such a scheme may be of value, this is not what the
author meant in the previous para.

\vspace{.12in} \noindent 4. The first three items in section 3
(relating to observables, states and the CC condition) were planned
to constitute a reasonably standardized noncommutative probabilistic
setting which, as we have seen, holds promise for being the proper
replacement of the deterministic setting of classical mechanics for
the description of dynamics of systems and, more generally, for
probability theoretic developments.

\vspace{.12in} \noindent 5. In the `mixed' case, when one of the
algebras, say $\sca^{(1)}$, is commutative and the other
noncommutative, it is possible to have a generalized symplectic
structure (of the type mentioned in section 1.5). Writing $fA$ for
$f \otimes A$, a general element of the tensor product algebra \sca
\ is of the form $\sum f_i A_i$ (finite sum); the product in \sca \
takes the form \begin{eqnarray*}  (\sum_i f_i A_i) ( \sum_j g_j B_j)
= \sum_{i,j} f_i g_j A_i B_j. \end{eqnarray*} The subalgebra
$\tilde{\sca}^{(1)}$ belongs to the center of \sca. Taking, in the
notation of section 1.5, $ \scx = IDer(\sca)$, we can have the
generalized symplectic algebra $(\sca, \scx, \omega)$ with $\omega =
b \omega_c$ giving the PB
\begin{eqnarray} \{ fA, gB \} = b^{-1} fg [A,B]. \end{eqnarray}
When $ \sca^{(1)}$ represents a classical system and $\sca^{(2)}$ a
quantum one, such an approach clearly amounts to treating the
classical observables as external fields. This is not adequate for a
proper treatment of the interaction of a classical and a quantum
system.

\vspace{.6in} \noindent \begin{tabbing} \hspace{.2in} \=
\scriptsize{INDIAN STATISTICAL INSTITUTE, SJS  SANSANWAL  MARG, NEW
DELHI 110016} \\
\>\small{\emph{E-mail address}: tulsi@iitk.ac.in}
\end{tabbing}

\newpage \begin{center} REFERENCES \end{center}

\vspace{.12in} \noindent
\begin{description}
\item[[1]] L. ACCARDI, Topics in quantum probability, \emph{Phys.
Rep.} \textbf{77} (1981), 169-192.
\item[[2]] L.M. ALONSO, Group-theoretic foundations of classical
and quantum mechanics. II. Elementary systems, \emph{J. Math. Phys.}
\textbf{20} (1979), 219-230.
\item[[3]] M.BORN, Zur quantenmechanik der stossvorg$\ddot{a}$nge,
\emph{Zeit. f. Physik} \textbf{37} (1926), 863-867.
\item[[4]] M. BORN and P. JORDAN, Zur quantenmechanik, \emph{Zs. f. Phys.}
\textbf{34} (1925), 858-888.
\item[[5]] M. BORN, W. HEISENBERG and P. JORDAN, Zur quantenmechanik II,
\emph{Zs. f.Phys.} \textbf{35}, 557-615.
\item[[6]] J.CARO and L.L. SALCEDO, Impediments to mixing classical and
quantum dynamics, \emph{Phys. Rev.} \textbf{A 60}, 842-852. [arXiv :
quant-ph/9812046.]
\item[[7]] A. CONNES, \emph{Noncommutative Geometry}, Academic
Press, New York, 1994.
\item[[8]] TULSI DASS, Noncommutative geometry and unified formalism
for classical and quantum mechanics, Indian Institute of Technology,
Kanpur preprint, 1993.
\item[[9]] TULSI DASS, Towards an autonomous formalism for quantum
mechanics, arXiv : quant-ph/0207104.
\item[[10]] TULSI DASS, Consistent quantum-classical interaction and
solution of the measurement problem in quantum mechanics, arXiv :
quant-ph/0612224.
\item[[11]] TULSI DASS, Supmech: a unified symplectic view of
physics, to be published.
\item[[12]] P.A.M. DIRAC, The fundamental equations of quantum
mechanics, \emph{Proc. Roy. Soc.} \textbf{A 109} (1926), 642-653.
\item[[13]] P.A.M. DIRAC, \emph{Principles of Quantum Mechanics},
Oxford University Press, London (1958).
\item[[14]] A.E.F. DJEMAI, Introduction to Dubois-Violette's
noncommutative differential geometry, \emph{Int. J. Theor. Phys.}
\textbf{34} (1995), 801-887.
\item[[15]] D.A. DUBIN and M.A. HENNINGS, \emph{Quantum Mechanics,
Algebras and Distributions} Longman Scientific and Technical,
Harlow, 1990.
\item[[16]] M. DUBOIS-VIOLETTE, Noncommutative differential
geometry, quantum mechanics and gauge theory, in \emph{Lecture Notes
in Physics, vol. 375}, Springer, Berlin, 1991, 13-24.
\item[[17]] M. DUBOIS-VIOLETTE, Some aspects of noncommutative
differential geometry, arXiv : q-alg/9511027.
\item[[18]] M. DUBOIS-VIOLETTE, Lectures on graded differential
algebras and noncommutative geometry, arXiv : math.QA/9912017.
\item[[19]] M. DUBOIS-VIOLETTE, R. KERNER and J. MADORE,
Noncommutative differential geometry of matrix algebras, \emph{J.
Math. Phys.} \textbf{31}, 316-322.
\item[[20]] I.M. GELFAND and N.J. VILENKIN, \emph{Generalized
Functions, vol. 4}, Academic press, New York, 1964.
\item[[21]] W. GREUB, \emph{Multilinear Algebra, 2nd edition},
Springer Verlag, New York, 1978.
\item[[22]] W. HEISENBERG, $\ddot{U}$ber quantentheoretische umdeutung
kinematischer und mechanischer beziehungen, \emph{Zs. f. Phys.}
\textbf{33}, 879-893.
\item[[23]]D. HILBERT, ``Mathematical problems'', lectures delivered
before the International Congress of Mathematicians in Paris in
1900, translated by M.W. NEWSON, \emph{Bull. Amer. Math. Soc.}
\textbf{8} (1902), 437.
\item[[24]] A. INOUE, \emph{Tomita-Takesaki Theory in Algebras of
Unbounded Operators}, Springer, Berlin, 1998.
\item[[25]] J.M. JAUCH, \emph{Foundations of Quantum Mechanics},
Addison-Wesley, Madison, Mass., 1968.
\item[[26]] P. JORDAN, J. VON NEUMANN and E. WIGNER, On an algebraic
generalization of the quantum mechanical formalism, \emph{Ann.
Math.} \textbf{35} (1934), 29-64.
\item[[27]] G.W. MACKEY, \emph{Mathematical Foundations of Quantum
Mechanics}, Benjamin-Cummings, Reading, Mass. 1968.
\item[[28]] P -A MEYER, \emph{Quantum Probability for Probabilists},
second edition, Springer, Berlin, 1995.
\item[[29]] J.E. MOYAL, Quantum mechanics as a statistical theory,
\emph{Proc. Camb. Phil. Soc.} \textbf{45} (1949), 99-124.
\item[[30]] J.  VON NEUMANN, \emph{Mathematical Foundations of
Quantum Mechanics}, Princeton University Press, 1955.
\item[[31]] K.R. PARTHASARATHY, \emph{An Introduction to Quantum
Stochastic Calculus}, Birkha$\ddot{u}$ser, Basel, 1992.
\item[[32]] E. SCHR$\ddot{O}$DINGER, Quantisierung als
eigenwertproblem, \emph{Annalen der Physik} \textbf{79} (1926),
361-376.
\item[[33]] I.E. SEGAL, Postulates for general quantum mechanics,
\emph{Ann. of Math.} \textbf{48} (1947), 930-948.
\item[[34]] E.C.G. SUDARSHAN and N. MUKUNDA , \emph{Classical
Dynamics : a Modern Perspective}, Wiley, New York, 1974.
\item[[35]] V.S. VARADARAJAN, \emph{Geometry of Quantum Theory},
second edition, Springer-Velag, New York, 1985.
\item[[36]] C.A. WEIBEL, \emph{An introduction to Homological
Algebra}, Cambridge University Press, 1994.
\item[[37]] H. WEYL, \emph{Theory of Groups and Quantum Mechanics},
Dover, New York, 1949.
\item[[38]] J.A. WHEELER and W.H. ZUREK, \emph{Quantum Theory of
Measurement}, Princeton University Press, 1983.
\item[[39]] E.P. WIGNER, On the quantum correction for thermodynamic
equilibrium, \emph{Phys. Rev.} \emph{40} (1932), 749-759.
\item[[40]] E.P. WIGNER, On unitary representations of the
inhomogeneous Lorentz group, \emph{Ann. of Math.} \textbf{40}
(1939), 149-204.
\item[[41]] N. WOODHOUSE, \emph{Geometric Quantization}, Clarendon
Press, Oxford, 1980.

\end{description}

\end{document}